# Design of Carbon Capture Processes Under Part-load Operating Conditions

**David Y. Shu[ab], Boxun Huang[a], Yurim Kim[a], Randall Field[b], Rahul Gandhi[c], and Sungho Shin[a*]**

[a] Massachusetts Institute of Technology, Department of Chemical Engineering, Cambridge, MA 02139, USA
[b] Massachusetts Institute of Technology, MIT Energy Initiative, Cambridge, MA 02139, USA
[c] Shell Technology Center, Houston, TX 77082, USA
* Corresponding Author: sushin@mit.edu

## ABSTRACT

Solvent-based carbon capture can reduce $CO_2$ emissions resulting from a continued reliance on fossil power plants for firm power. These capture processes remove $CO_2$ from flue gases via a solvent. Careful design via process systems optimization can limit the overall cost of carbon capture, which is both capital- and energy intensive. As dispatchable power plants operate to meet varying load demand, the design process needs to account for varying operating points. However, optimizing the design over multiple operating points yields high computational complexity, which is why designs are often based on a single operating point in practice. Here, we identify optimal carbon capture process designs via stochastic optimization, reducing computational complexity through a data-driven approach-to-equilibrium model of the absorption and desorption processes. We represent variable flue gas conditions based on part-load operation data of a representative coal power plant. Accounting for this variability in the design substantially reduces equipment size and total plant cost by 6–9 % at the expense higher operating costs, yielding a reduction in total cost of carbon capture by 0.7–1.7 %. Given the capital intensity of carbon capture, variability of flue gas conditions therefore should be considered at the design stage, particularly if capture is deployed on plants subject to load following.

## INTRODUCTION

Fossil power plants accounted for almost 60 % of power generation in the US in 2025 [1]. Renewed interest in natural gas power plants to supply firm power to data centers may lock in additional $CO_2$-emissions over their lifetime. These emissions can be reduced by carbon capture processes.

Carbon capture processes are both capital- and energy-intensive: Mature capture processes rely on solvents, removing $CO_2$ from flue gases through contacting with a solvent in countercurrent flow in absorber columns (Fig. 1). The $CO_2$-rich solvent is regenerated in a desorber using steam, producing a $CO_2$ stream that needs to be conditioned for permanent sequestration or utilization.

Computational process system engineering is a crucial tool for designing capture processes that balance capital and operational costs, while meeting process targets and technical and regulatory constraints. Such designs should account for uncertainties which are prevalent in operational conditions and in economic and technical parameters [2]. In capture processes, uncertainties arise, e.g., due to the varying flue gas characteristics of power plants: As fossil power plants provide grid-balancing flexibility, they are increasingly operated at part load with growing renewable penetration. Accordingly, carbon capture plants should be designed for a wide range of part-load conditions.

However, accounting for uncertainty in process systems design can be computationally challenging [2]. To avoid excessive complexity in carbon capture design, rigorous handling of uncertainties is sometimes avoided by applying overdesign margins to designs obtained from deterministic simulation [3] or from optimization with detailed parametric studies [4]. Alternatively, computational complexity can be reduced by applying surrogate models [5,6], or by focusing on select components [7,8].

Uncertainty can be handled via robust optimization, yielding designs that guarantee feasibility for considered uncertainty sets [7,9,10] but may be overly conservative [2]. Alternatively, uncertainty can be addressed by stochastic programming which is risk-neutral [2].

Prior research already considered the optimal sizing of the desorber and storage tanks via stochastic programming considering uncertain electricity and $CO_2$ prices or taxes [8]. Seo et al. [11] further applied stochastic optimization to the design of carbon capture using typical days as scenarios, focusing on dynamic operation to respond to a varying electricity price and demand.

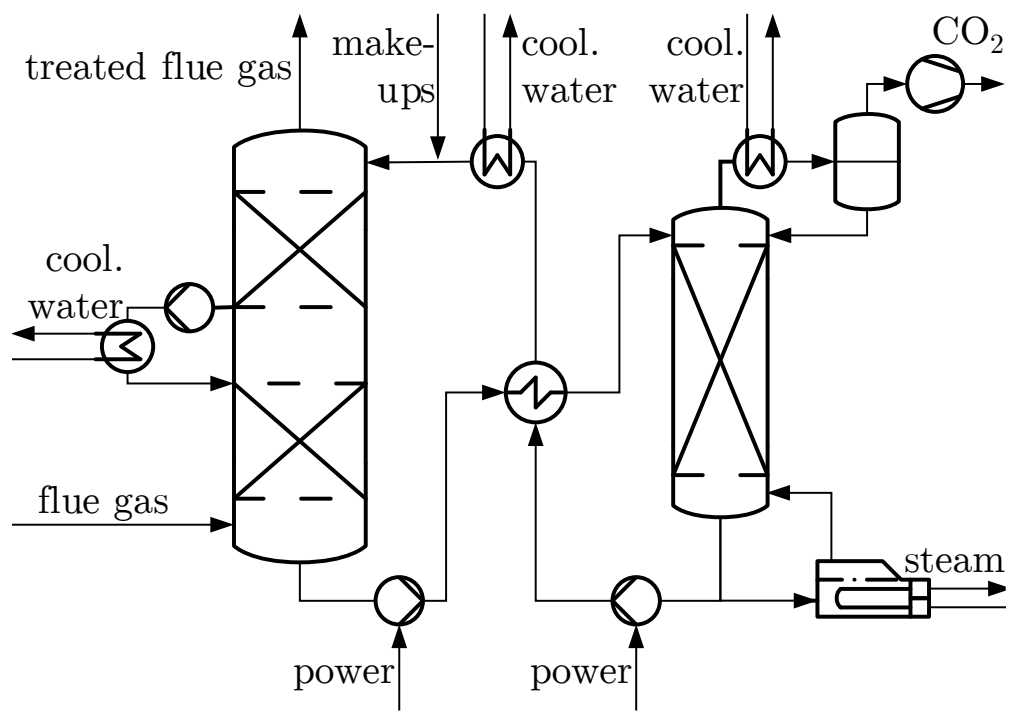

**Figure 1.** Process flow diagram of the carbon capture process considered for design optimization.

In contrast to previous studies, we apply two-stage stochastic optimization to isolate the impact of part-load flue gas variability as a design driver for capture plants. To this end, we derive scenarios from real-world distributions of part-load operation levels of power plants. We model the capture process as a set of interconnected unit operations, and include detailed technoeconomic models, allowing us to minimize the total annualized cost, while achieving capture targets.

# MODELING AND METHODOLOGY

We formulate the design of the carbon capture process as a stochastic programming problem (((1)–((4)), assuming 30 wt% monoethanolamine (MEA) as solvent.

$$\min_{\boldsymbol{x},\{\boldsymbol{y}_s\}_{s\in\mathcal{S}}} c^{\mathrm{cap}}(\boldsymbol{x}) + \sum_{s\in\mathcal{S}} w_s \; c^{\mathrm{op}}(\boldsymbol{x},\boldsymbol{y}_s) \tag{1}$$

$$\text{s.t.} \quad \boldsymbol{f}_s(\boldsymbol{x},\boldsymbol{y}_s) \leq 0, \;\; s \in \mathcal{S} \tag{2}$$

$$\boldsymbol{h}(\boldsymbol{x},\{\boldsymbol{y}_s\}_{s\in\mathcal{S}}) \leq 0, \tag{3}$$

$$\boldsymbol{l_x} \leq \boldsymbol{x} \leq \boldsymbol{u_x}, \quad \boldsymbol{l_y} \leq \boldsymbol{y}_s \leq \boldsymbol{u_y}, \quad s \in \mathcal{S} \tag{4}$$

The objective function ((1)) is the total annualized costs, where $c^{\mathrm{cap}}(\boldsymbol{x})$ represents the annualized total ownership cost of the equipment, $c^{\mathrm{op}}(\boldsymbol{x},\boldsymbol{y}_s)$ the scenario-dependent operating costs, and $w_s$ the scenario weight. $\boldsymbol{f}_s$ and $\boldsymbol{h}$ are scenario-dependent and scenario-coupling constraints, respectively, representing, e.g., techno-economic models, capture rate targets, and technical process models, such as thermodynamic properties, unit process models, connectivity of the unit processes, and balance equations. $\boldsymbol{l}_x$, $\boldsymbol{u}_x$ and $\boldsymbol{l}_y$, $\boldsymbol{u}_y$ are the bounds of the first and second stage variables.

The first-stage design variables $\boldsymbol{x}$ size the unit-operations in the capture process (Fig. 1), specifically the packing height and diameter of the columns, maximum pump flow rates, number of channels in the cross-heat exchanger, maximum compressor power, and heat-exchanger surface areas of reboiler, condenser, cooler, and intercooler. We relax the number of channels in the cross heat exchanger to a continuous variable to maintain a nonlinear programming formulation. The second-stage operating decisions $\boldsymbol{y}_s$ for scenarios $s \in \mathcal{S}$, include the operational variables, e.g., reboiler duty, solvent loading, or solvent flow rate.

## Technoeconomic Model

The bare erected costs are estimated based on techno-economic models [12], using the component sizes from the process model and assuming equipment made of stainless steel. The column heights are estimated using the same column to packing height ratio reported in [13]. The diameter is estimated based on a generalized pressure drop correlation for structured packing [14] at an 80 % approach to flooding. The costs of structured packing material are estimated assuming *Mellapak 250Y* [14]. The cross heat exchanger is costed as a plate-and-frame heat exchanger, while reboiler, condenser, cooler, and intercooler are estimated based on tube-and-shell heat exchangers. As the required size of the heat exchangers and variable speed pumps exceed the maximum sizes specified in [12], we assume multiple maximum-size units in parallel. For those units, the cost is estimated by scaling the cost of one maximum-size unit linearly with the number of units needed. We relax the number of units needed to a continuous variable to maintain a nonlinear programming formulation. The compressor power and cooling loads are estimated based on [15].

The bare erected costs are escalated to 2023 values using the Chemical Engineering Plant Cost Index. We estimate engineering, procurement and construction costs and contingencies to arrive at total plant costs and further include preproduction, inventory, financing, and other costs using factors based on [16,17] to arrive at total ownership cost. The total ownership cost is annualized using the capital charge factor from [17].

The operating costs include the utility costs of steam at 150 psig, electricity, and cooling water at 90 °F [12]. We further include fixed operating costs for labor, maintenance, taxes and insurance, as well as variable operating costs for waste treatment and chemicals, scaled by the amount of CO2 captured [16]. The cost of MEA makeup is estimated based on degradation rates in [21] and prices in [22] escalated to 2023 [23].

## Process Model

The component sizing, utility demands, and amount of $CO_2$ captured is determined using the process system model. The process models include material and energy balances. As a simplifying assumption, we calculate enthalpy changes in liquid and vapor phase using the pure-component models for heat capacity, heat of vaporization and heat of absorption [20], neglecting heat of mixing. In addition, we include equilibrium reactions of MEA with $CO_2$ [24] in the absorber and desorber columns, condenser, reboiler, and cross heat exchanger.

As a simplifying approximation, we assume vapor-liquid equilibrium in condenser and reboiler based on Henry's law for $CO_2$ [25] and Raoult's law for water using the vapor pressure [26]. For simplicity, we neglect the solubility of $O_2$ and $N_2$ and assume MEA to be nonvolatile.

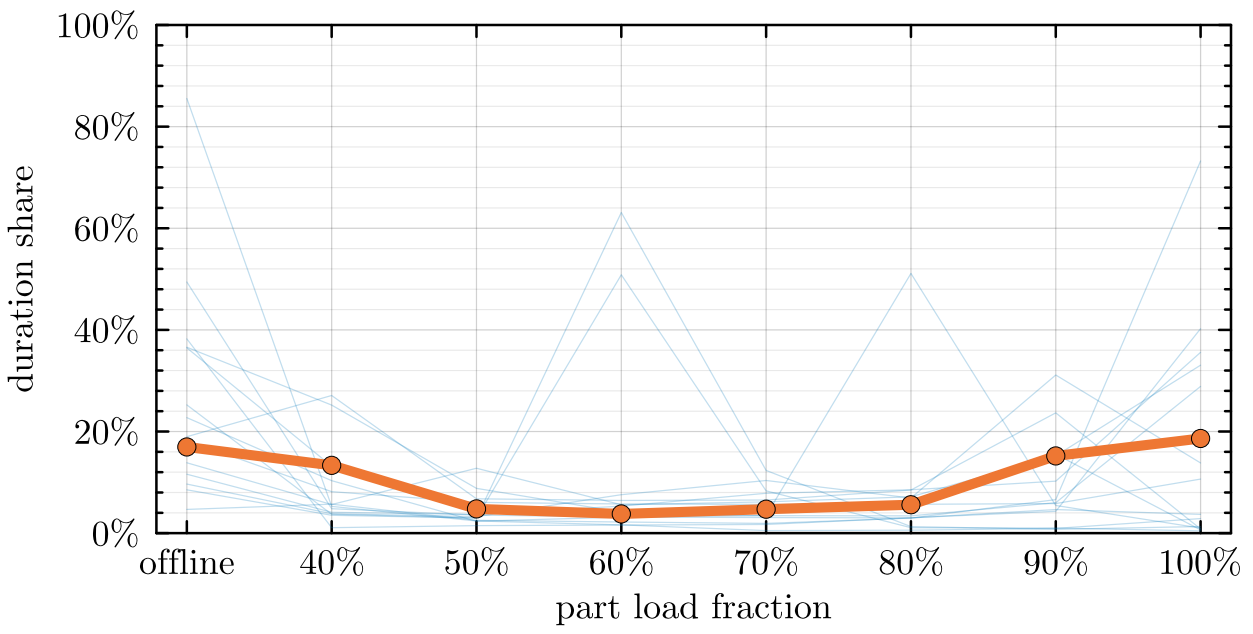

**Figure 2.** Duration in part-load operation of sub-bituminous coal power plant units in Texas in 2025 [18]. A representative operation profile determined via k-medoids clustering is highlighted in orange.

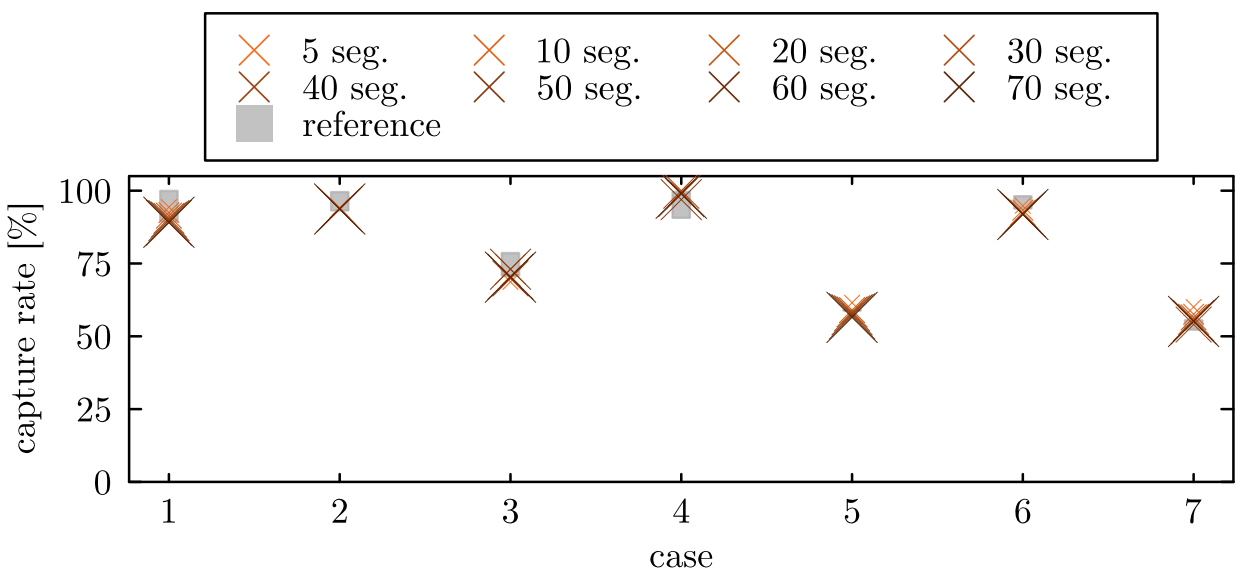

**Figure 3.** Simulated capture rates in comparison to measurements (reference) for validation. The inlet flue gas and lean solvent conditions of the measurements are summarized by [19,20]. Simulation results are shown for a varying number of column segments (seg.).

To account for rate-limitations of chemical reactions and mass transfer between the liquid and vapor phase in the columns, we use an approach-to-equilibrium model (((5)) to determine how much $CO_2$ is absorbed or desorbed based on the equilibrium driving force. This modeling approach reduces problem complexity in comparison to a more rigorous rate-based model. We discretize the columns axially into segments $j$, which describe local vapor-liquid mass transfer and chemistry. Segments are numbered from the bottom of the column, and the approach-to-equilibrium model is applied to each $j$:

$$y_{\mathrm{CO_2},s,j+1} = y_{\mathrm{CO_2},s,j} - \eta(k_{s,j}, \Delta z)\left(y_{\mathrm{CO_2},s,j} - y^*_{\mathrm{CO_2},s,j}\right) \quad (5)$$

Here, $y^*_{\mathrm{CO_2},s,j}$ is the $CO_2$ mole fraction in the vapor phase at equilibrium and the function η depends on parameter $k_{s,j}$ and segment height $\Delta z$. We choose the function form so that η remains in $[0,1]$ and is consistent across discretization resolutions (((6)).

$$\eta(k, \Delta z) = 1 - e^{-k\Delta z} \quad (6)$$

Parameter $k_{s,j}$ lumps the effects of gas- and liquid-phase mass transfer resistance, chemical reaction rates, interfacial area and hydrodynamics, and modeling approximations on $CO_2$ absorption and desorption. Nevertheless, the liquid phase composition is always based on reaction equilibrium. The parameter $k_{s,j}$ is obtained by evaluating the calibrated function $k(\cdot)$ at the local conditions $\boldsymbol{\theta}_{s,j}$: $k_{s,j} = k(\boldsymbol{\theta}_{s,j})$, where $\boldsymbol{\theta}_{s,j}$ includes the vapor-to-liquid molar flow rate ratio, inlet temperature, pressure, inlet solvent $CO_2$ loading, MEA mass fraction in the $CO_2$-free inlet solvent, and superficial gas velocity. The calibrated function $k(\cdot)$ has the form:

$$\ln k = \beta_0 + \sum_{n\in\mathcal{N}} (\beta_n \theta_n) + \sum_{(n,m)\in\mathcal{M}} (\beta_{n,m}\, \theta_n\, \theta_m) \quad (7)$$

Here, $\mathcal{N}$ is the index set of calibration variables and $\mathcal{M} \subseteq \mathcal{N} \times \mathcal{N}$ is a set of index pairs included to capture interactions. The coefficients $\beta_0$, $\beta_n$, and $\beta_{n,m}$ are fitted to calibration data from experiments or detailed simulations. For clarity, (7) omits the normalization of the calibration data using mean and standard deviation.

The plate-and-frame cross heat-exchanger is modeled following [27]. The tube-and-shell heat exchangers are sized assuming constant overall heat transfer coefficient from [28], and using Chen's approximation [29] of the logarithmic mean temperature difference. Finally, we assume an intercooler in the center of the absorber column with a target temperature range of 40 °C–50 °C. For density calculations, we assume ideal gases and constant solvent densities of 1000 kg $m^{-3}$ for simplicity.

The capture rate targets is specified for every operating point ((8)) or as an average over all scenarios ((9)).

$$V_s^{\mathrm{feed}}\, y_{\mathrm{CO_2},s}^{\mathrm{feed}} - V_s^{\mathrm{abs,out}}\, y_{\mathrm{CO_2},s}^{\mathrm{out}} \geq \gamma^{\mathrm{cap}} V_s^{\mathrm{feed}}\, y_{\mathrm{CO_2},s}^{\mathrm{feed}}, \;\forall s \in \mathcal{S} \quad (8)$$

$$\sum_{s\in\mathcal{S}} w_s\left(V_s^{\mathrm{feed}}\, y_{\mathrm{CO_2},s}^{\mathrm{feed}} - V_s^{\mathrm{abs,out}}\, y_{\mathrm{CO_2},s}^{\mathrm{out}}\right) \geq \gamma^{\mathrm{cap}} \sum_{s\in\mathcal{S}} w_s\left(V_s^{\mathrm{feed}}\, y_{\mathrm{CO_2},s}^{\mathrm{feed}}\right) \quad (9)$$

Here, $V_s^{\mathrm{feed}}$ and $V_s^{\mathrm{abs,out}}$ are the absorber vapor inlet and outlet molar flows, and $y_{\mathrm{CO_2},s}^{\mathrm{feed}}$ and $y_{\mathrm{CO_2},s}^{\mathrm{out}}$ the $CO_2$ molar fractions. $\gamma^{\mathrm{cap}}$ is the capture rate target.

# CASE STUDY

In the case study, we solve the design and operation problem (((1)–(4)) of the carbon capture process (Fig. 1), considering flue gases from coal power plants. For simplicity, we assume an end-of-pipe process and no changes to power plant operation or performance.

To quantify the value of stochastic design, we compare four cases: As a reference, we design the process assuming a single scenario at the design operating point with 80 % power plant utilization (*deterministic*), reflecting design practice in industry. The operating costs are evaluated by optimizing the operation for the fixed deterministic design.

**Table 1.** Optimal sizes of the main process equipment in the deterministic cases and relative changes in the stochastic cases.

| variable | deterministic reference | stochastic | stochastic (avg.) |
|---|---|---|---|
| packed height (abs.) | 27.2 m | -10.8 % | -11.3 % |
| diameter (abs.) | 11.8 m | +1.0 % | -0.5 % |
| packed height (des.) | 6.1 m | -13.6 % | -12.3 % |
| diameter (des.) | 7.9 m | +2.4 % | -1.1 % |
| cross HEX area | 4.4×10$^4$ m$^2$ | -46.8 % | -42.6 % |
| condenser area | 7.3×10$^4$ m$^2$ | -3.8 % | -8.6 % |
| cooler area | 3.9×10$^4$ m$^2$ | +17.3 % | +3.7 % |
| intercooler area | 6.2×10$^4$ m$^2$ | -8.8 % | -12.9 % |
| reboiler area | 4.4×10$^4$ m$^2$ | +4.4 % | -1.5 % |
| compressor power | 17.8 MW | 0.0 % | -3.0 % |
| lean pump flow | 1.8×10$^4$ gpm | +5.3 % | -3.0 % |
| rich pump flow | 1.9×10$^4$ gpm | +4.9 % | -3.0 % |

Next, we calculate a stochastic design, considering multiple weighted load scenarios (*stochastic*). In both, the *stochastic* and the *deterministic* cases, a 95 % capture rate needs to be achieved in all scenarios ((8)). We further include deterministic and stochastic cases in which we assume that the capture rate needs to be achieved on average across all scenarios (*deterministic (avg.)* and *stochastic (avg.),* (9)).

The scenarios $s \in \mathcal{S}$ represent operation at power plant loads from 40–100 % at increments of 10 %, and differ in the flue gas characteristics: We include flue gas flows and $CO_2$ concentration at various operating loads for a representative coal power plant [30]. Due to the large flue gas mass flow of the plant, we assume two parallel carbon capture process trains. For simplicity, we assume the flue gas stream to be evenly distributed between trains, neglecting the option to shut down individual trains at low load or shared equipment between trains. As the $H_2O$ and $O_2$ concentrations of the flue gas are not reported, we use constant values from [16].

The relative duration at each operating load is based on a representative coal power plant. Specifically, we analyzed the operating points of sub-bituminous coal power plants in Texas in 2025 [18]: For each power plant, we allocate all operating data points to the nearest 10 %-part-load bin, to estimate the duration share in each part-load operating point. We then determine the medoid to select a power plant with representative operation profile for the scenario weights in the case study (Fig. 2).

The calibration coefficients $\beta_0$, $\beta_n$, and $\beta_{n,m}$ ((7)) are obtained by a least-squares fit to 190,000 datapoints simulated using the *CCSI steady state MEA model* [31], by sampling from a range of operating conditions relevant for the case study. Predicting the $CO_2$ concentration at the segment output using the calibrated model results in an average prediction error of 1.5 %. The 5th–95th percentile prediction error ranges from -3.0 to 6.5 %.

We validate the carbon capture process model by simulating a pilot plant capture process from literature [19,20]. For validation, we fix the flue gas and lean solvent characteristics, and the capture plant design based on the cases reported in [20], and compare our simulated capture rates to the ones reported in [20]. The simulated capture rates show agreement with the measurements, with the number of segments specified in our model having only a minor impact on the results (Fig. 3).

To solve the design optimization, we obtain good initial values by simulating an open-loop capture process for a single scenario at a column resolution of 5 segments. To do so, we tear the lean-solvent recycle stream, and fix the lean-solvent characteristics, the plant design, and the reboiler duty at reasonable values. The simulation results then serve as starting values for the optimization problem of the same resolution. From there, the number of scenarios and the column resolution are increased successively, using the previous solution as initial values and applying linear interpolation when increasing the column resolution. We select 30 segments as the final resolution for both columns, which is sufficiently high for the convergence of the capture rate (Fig. 3).

The optimization problems are implemented in *Julia* using *JuMP* v1.29.3 and *IPOPT* v1.13.0 with the linear solver *MA57*. Thus, our results report local optima.

# RESULTS

The optimal size of most units is reduced in the two stochastic cases compared to the deterministic case (Tab. 1): Absorber and desorber packing heights are lowered by 11 – 13 %, and cross heat exchanger area is reduced by up to 47 %. Smaller equipment reduces total plant cost by 6–9 % (Fig. 4 (a)). This reduction is driven by the columns and heat exchangers, which make up65 % of the cost in the deterministic case.

The reductions in total plant cost are partially offset by higher operating costs. The operating costs are dominated by regeneration steam consumption, which increases in the stochastic cases by 4–5 % (Fig. 4 (b)) due to the smaller equipment, particularly the cross heat exchanger. Overall, the stochastic cases only reduce total annualized cost by 0.7–1.7 % (Fig. 4 (c)).

Reductions of total annualized cost are slightly higher in the case in which the 95 % capture rate is achieved on an annual average. Here, additional operational flexibility is exploited by shifting $CO_2$ capture from full-load to part-load periods: In the stochastic case, the capture rate increases from 92 % at full load to up to 98 % at part load. This shift reduces peak mass flows and lowers total plant cost in the stochastic case. The deterministic case does not benefit from the additional flexibility, as its design is based on the design operating point.

# CONCLUSIONS

This study applies two-stage stochastic optimization to simultaneously optimize the design and operation of a solvent-based carbon capture process system under variable flue gas conditions. We model absorption and desorption

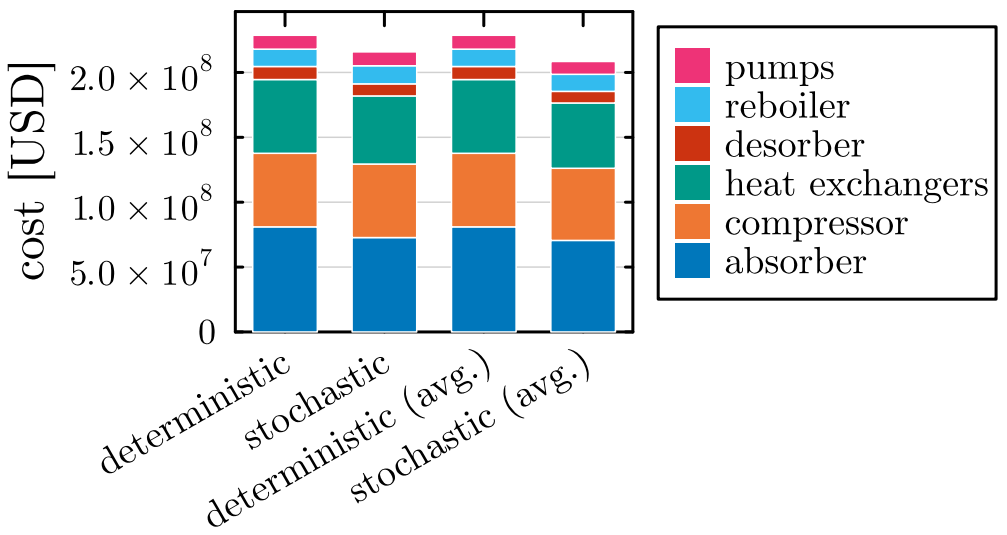

(a) Total plant cost of the capture system.

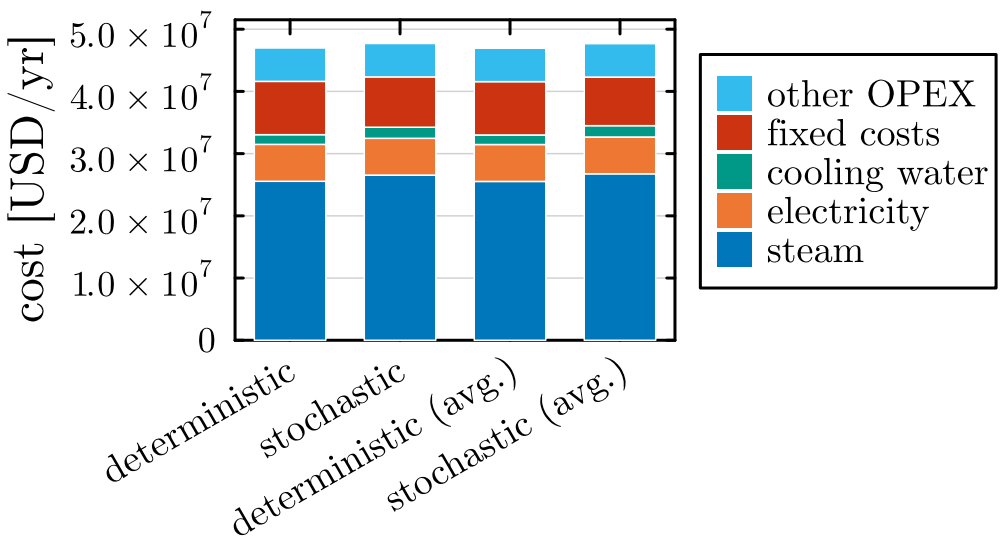

(b) Annual operating cost of the capture system.

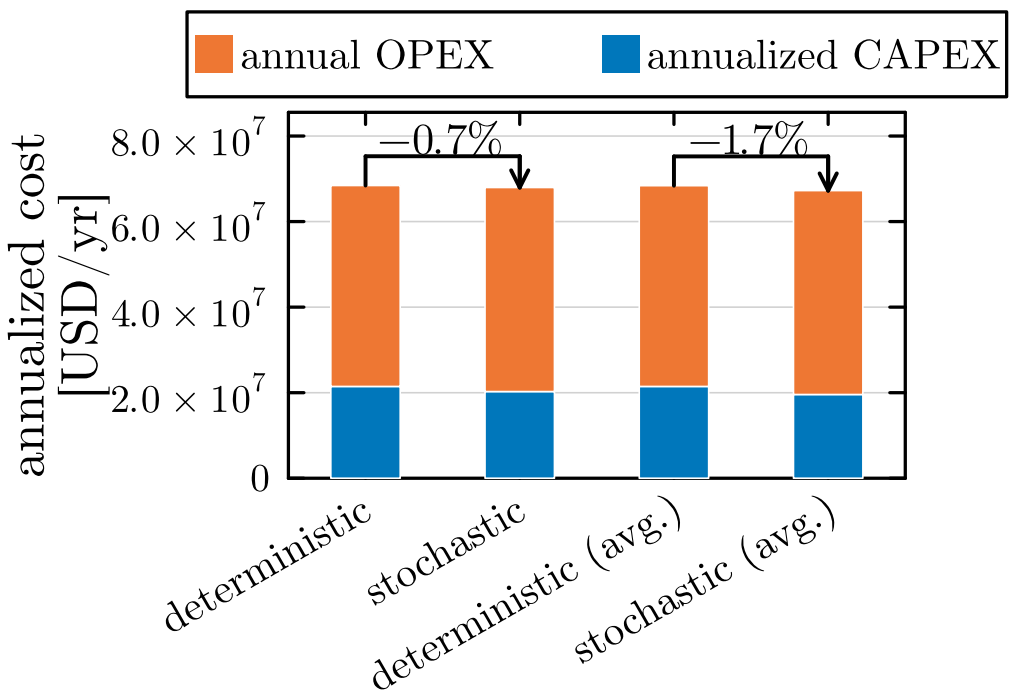

(c) Total annualized cost of the capture system.

**Figure 4.** Cost overview of total plant cost and operating costs for the deterministic and stochastic designs, assuming 95 % capture rate in every scenario or on average over the full year (avg.).

processes in the columns using a data-driven approach-to-equilibrium model to approximate rate-limited capture. The total annualized cost is calculated using technoeconomic models as the objective function.

Our case study of carbon capture from a coal power plant includes scenarios differing in flue gas conditions with scenario weights derived from the part-load operation of a representative power plant in Texas in 2025. The stochastic approach accounts for multiple part-load operating points and results in smaller sizes for most equipment, which translates into substantial reductions of the total plant cost by 6–9 %. Total annualized costs are reduced by 0.7–1.7 % compared to the deterministic approach. Targeting an annual capture rate introduces additional flexibility that allows to shift $CO_2$ capture from full-load to part-load periods, enabling a further reduction in equipment size and total plant cost.

This study demonstrates the cost reduction potential of stochastic optimization in carbon capture process design. In future research, the stochastic problem formulation can be extended by parametric uncertainty in plant performance. Further, the capture process can be co-designed with the power plant, enabling cost reductions through integration. Crucially, the impacts of steam extraction can be explored through co-design. The increase in computational complexity resulting from these extensions needs to be addressed, e.g., via GPU-accelerated optimization [32].

## ACKNOWLEDGEMENTS

DS, BH, and SS received support from the Future Energy Systems Center through the MIT Energy Initiative. This publication is based on work supported by the U.S. Department of Energy under contract DE-AC02-06CH11357. Funding was provided by the U.S. DOE Hydrocarbons and Geothermal Energy Office under FWP-DOE-145-SE-14205HPC4E-25.

## DECLARATION OF USE OF AI

The authors used artificial intelligence tools to assist with grammar, spelling, and readability improvements in the manuscript. No AI tools were used to generate scientific content, analysis, or conclusions. All technical material reflects the authors' original work, and the authors are fully responsible for the accuracy of the content.